\documentclass[11pt]{article}

\usepackage[compact]{titlesec}
\usepackage{array}
\usepackage{amsfonts}
\usepackage{amscd}
\usepackage{amssymb}
\usepackage{amsthm}
\usepackage{amsmath}
\usepackage{stmaryrd}
\usepackage{graphicx}
\usepackage{color}
\usepackage{verbatim}
\usepackage{mathrsfs}
\usepackage{pstricks,tikz}
\usetikzlibrary{calc}
\usepackage{caption}
\usepackage[neveradjust]{paralist}
\usepackage[a4paper, total={15cm, 25cm},centering]{geometry}
\usepackage{float}
\usepackage[absolute]{textpos}
\usepackage{subcaption}
\usepackage{ytableau}
\usepackage{hyperref}
\usepackage{bbm}
\usepackage{lineno}

\hypersetup{hidelinks}

\usetikzlibrary{automata, positioning, arrows}

\let\OLDthebibliography\thebibliography
\renewcommand\thebibliography[1]{
  \OLDthebibliography{#1}
  \setlength{\parskip}{0pt}
  \setlength{\itemsep}{0pt plus 0.3ex}
}

\definecolor{green1}{RGB}{153,216,201}
\definecolor{green2}{RGB}{44,190,95}
\definecolor{ForestGreen}{cmyk}{0.91,0,0.88,0.42}

\definecolor{blue1}{RGB}{158,202,225}
\definecolor{blue2}{RGB}{49,130,189}

\definecolor{RedOrange}{cmyk}{0,0.77,0.87,0} 
\definecolor{Mahogany}{cmyk}{0,0.85,0.87,0.35} 
\definecolor{Maroon}{cmyk}{0,0.87,0.68,0.32} 
\definecolor{BrickRed}{cmyk}{0,0.89,0.94,0.28} 
\definecolor{Red}{cmyk}{0,1.,1.,0} 
\definecolor{OrangeRed}{cmyk}{0,1.,0.50,0} 

\definecolor{purple}{rgb}{0.8,0.12,0.8}
\definecolor{orange}{rgb}{1.0,0.7,0.0}
\definecolor{pink}{rgb}{1,0.5,0.8}
\definecolor{blackg}{rgb}{0.1,0.25,0.1}
\definecolor{ForestGreen}{cmyk}{0.91,0,0.88,0.42}
\definecolor{Turquoise}{cmyk}{0.85,0,0.20,0}

 \theoremstyle{plain}

\newtheorem{thm}{Theorem}[section]
\newtheorem{lemma}[thm]{Lemma}
\newtheorem{prop}[thm]{Proposition}
\newtheorem{cor}[thm]{Corollary}

\theoremstyle{definition}

\newtheorem{defn}[thm]{Definition}

\newtheorem{remark}[thm]{Remark}
\newtheorem{example}[thm]{Example}

\newtheorem{assumptions}[thm]{Assumptions}

\numberwithin{equation}{section}

\def\pr{\mathsf{pr}}
\def\tw{\widetilde{\mathsf{w}}_I}
\def\cact{\mathsf{Cact}}

\def\sA{\mathsf{A}}

\def\sC{\mathsf{C}}

\def\sF{\mathsf{F}}
\def\sG{\mathsf{G}}

\def\sT{\mathsf{T}}

\def\cD{\mathcal{D}}

\def\cJ{\mathcal{J}}

\def\cL{\mathcal{L}}

\def\cR{\mathcal{R}}

\newcommand{\la}{\lambda}

\newcommand{\Ga}{\Gamma}

\def\QQ{\mathbb{Q}}
\def\RR{\mathbb{R}}

\def\ZZ{\mathbb{Z}}

\DeclareMathOperator\Tr{\mathsf{Tr}}

\def\la{\lambda}

\def\<{\langle}
\def\>{\rangle}

\def\la{\lambda}

\def\sR{\mathsf{R}}

\newcommand{\sw}{\mathsf{w}}

\newcommand{\st}{\mathsf{t}}

\newcommand{\ba}{\mathbf{a}}

\newcommand{\sq}{\mathsf{q}}

\makeatletter
\renewcommand{\@makefnmark}{\mbox{\textsuperscript{}}}
\makeatother

\title{Cells, cacti, and the asymptotic algebra}
\author{F. Dos Santos and J. Parkinson}
\date{\today}

\begin{document}

\maketitle

\vspace{-0.5cm}

\begin{abstract}
We prove that the action of the cactus group on a Coxeter group defined by Losev and Bonnaf\'e can be realised as an action on the asymptotic algebra of the Hecke algebra, thus resolving a conjecture of Rouquier and White. We prove this in the setting of general (potentially unequal) parameters, conditional on Lusztig's conjectures~$\mathsf{P1}$--$\mathsf{P15}$, which are known to hold in the equal parameter case. 
\end{abstract}
%


\section{Introduction}

Let $(W,S)$ be a Coxeter system with $S$ finite, and let $\cact(W)$ be the associated cactus group. In~\cite{Los:19} Losev discovered a remarkable action of $\cact(W)$ on $W$ that is compatible with Kazhdan-Lusztig cells. More precisely, he defined an action of $\cact(W)\times\cact(W)$ on $W$ such that the ``left'' action preserves right cells and permutes left cells, and the ``right'' action preserves left cells and permutes right cells. Losev's work was in the setting of equal parameter finite Weyl groups, and in~\cite{Bon:16} Bonnaf\'e extended this action to arbitrary weighted Coxeter groups (under standard assumptions of a subset of Lusztig's properties $\mathsf{P1}$--$\mathsf{P15}$). 

Recently Rouquier and White~\cite{RW:24} showed that there is a natural linear action of $\cact(W)$ on the Hecke algebra (with extended coefficients), and conjectured that this action specialises to a $\ZZ$-linear action on Lusztig's asymptotic algebra~$\cJ$, and that this action recovers the original Losev-Bonnaf\'e action (see \cite[Conjecture~4.2]{RW:24}). In this paper we resolve the conjecture of Rouquier and White, thus providing a direct link between the cactus group and the asymptotic algebra. While Conjecture~4.2 of \cite{RW:24} is stated for equal parameter Hecke algebras, there is a natural extension to unequal parameters, and we will prove the conjecture in this generality under the assumption of properties $\mathsf{P1}$--$\mathsf{P15}$ (which are known to hold in the equal parameter case). 

Let us give a brief description of \cite[Conjecture~4.2]{RW:24}, and an indication of our approach. We first describe Bonnaf\'e's construction. Let $H=H(W,S)$ be a Hecke algebra defined over $\sR=\ZZ[\sq,\sq^{-1}]$ (possibly with unequal parameters). Let $(T_w\mid w\in W)$ and $(C_w\mid w\in W)$ be the standard basis, and the Kazhdan-Lusztig basis, of $H$ (respectively), leading to the notion of left, right, and two-sided cells in the (weighted) Coxeter group. In \cite{Bon:16} Bonnaf\'e constructs an action of $\cact(W)\times\cact(W)$ on $W$ such that if $\gamma^L$ (respectively~$\gamma^R$) denotes the permutation of $W$ obtained through the action of $(\gamma,1)\in\cact(W)\times\cact(W)$ (respectively $(1,\gamma)$) then:
\begin{compactenum}[$(1)$]
\item If $\Gamma$ is a left cell of $W$, then $\gamma^L(\Gamma)$ is also a left cell;
\item If $\Gamma$ is a right cell of $W$, then $\gamma^R(\Ga)$ is also a right cell;
\item If $w\in W$ then $\gamma^L(w)\sim_R w$ and $\gamma^R(w)\sim_L w$.
\end{compactenum}
Moreover, he shows that for each left (respectively right) cell $\Ga$, there exists a sign map such that the Kazhdan-Lusztig cell modules associated to $\Gamma$ and $\gamma^L(\Ga)$ (respectively $\gamma^R(\Ga)$) are isomorphic. 

The cactus group $\cact(W)$ is generated by elements $\gamma_I$ with $I\subseteq S$ a connected spherical subset, and so for each such $I$ Bonnaf\'e's construction gives permutations $\gamma_I^L$ and $\gamma_I^R$ of $W$, and associated sign maps $\eta_{I,L}:W\to\{-1,1\}$ and $\eta_{I,R}:W\to\{-1,1\}$.

The main theme of~\cite{RW:24} is to realise the Losev-Bonnaf\'e action as a ``crystal limit'' of a natural action of the cactus group on the Hecke algebra. It is well known that the map $\cact(W)\to W$ with $\gamma_I\mapsto \sw_I$ is a group homomorphism, where $\sw_I$ is the longest element of the parabolic subgroup~$W_I$. Now, writing $H_I$ for the spherical Hecke algebra $H(W_I,I)$, it is shown in \cite{RW:24} that the element $T_{\sw_I}\in H_I$ can be ``unitarised'' after extending scalars to produce a unique element $\tw\in \tilde{H}_I=\QQ(\sq)\otimes_{\sR}H_I$ satisfying various properties. Moreover, the properties of $\tw$ imply that homomorphism $\cact(W)\to W$ lifts to a homomorphism
$$
\theta:\cact(W)\to \tilde{H}^{\times}\quad\text{with}\quad \theta(\gamma_I)=\tw,
$$
where $\tilde{H}=\QQ(\sq)\otimes_{\sR}H$. Thus there is a natural $\QQ(\sq)$-linear action of $\cact(W)\times\cact(W)$ on $\tilde{H}$ given by 
$$
(\gamma,\gamma')\cdot h=\theta(\gamma)h\theta(\gamma')^{-1}\quad\text{for $\gamma,\gamma'\in\cact(W)$ and $h\in \tilde{H}$.}
$$

Let $\cJ$ be Lusztig's \textit{asymptotic algebra}. Thus $\cJ$ is a $\ZZ$-algebra with basis $(\st_w\mid w\in W)$, and multiplication given in terms of certain leading coefficients of the structure constants of the Kazhdan-Lusztig basis. Passing through Lusztig's homomorphism $\psi:H\to\ZZ[\sq,\sq^{-1}]\otimes_{\ZZ}\cJ$ (lifted to a homomorphism $\psi:\tilde{H}\to\QQ(\sq)\otimes_{\ZZ}\cJ$), the action of the cactus group on $\tilde{H}$ induces an action of the cactus group on $\tilde{\cJ}=\QQ(\sq)\otimes_{\ZZ}\cJ$. Specifically, the map
$$
\psi\circ\theta:\cact(W)\to \tilde{\cJ}^{\times}\quad\text{with}\quad \gamma_I\mapsto \psi(\tw)
$$
defines a group homomorphism, and there is a $\QQ(\sq)$-linear action of $\cact(W)\times\cact(W)$ on $\tilde{\cJ}$. 

Conjecture~4.2 of Rouquier and White~\cite{RW:24} concerns the possibility of a specialisation of the above action, producing an action of $\cact(W)\times\cact(W)$ on the $\ZZ$-algebra $\cJ$ recovering the original Losev-Bonnaf\'e action. More precisely, the conjecture (extended to the setting of unequal parameters) is as follows:
\smallskip

\begin{compactenum}[$(1)$]
\item The elements $\psi(\tw)\in\tilde{\cJ}$ admit specialisations at both $\sq=0$ and $\sq=1$ into~$\cJ$;
\item Writing $\psi_0(\tw)\in\cJ$ for the specialisation of $\psi(\tw)$ at $\sq=0$, for $w\in W$ we have 
$$
\psi_0(\tw)\st_w=\eta_{I,R}(w)\st_{\gamma_I^R(w)}\quad\text{and}\quad \st_w\psi_0(\tw)=\eta_{I,L}(w)\st_{\gamma_I^L(w)}. 
$$
\item We have
$$
\psi_0(\tw)=\sum_{d\in\cD}n_d\eta_{I,R}(d)\st_{\gamma_I^R(d)}=\sum_{d\in\cD}n_d\eta_{I,L}(d)\st_{\gamma_I^L(d)},
$$
where $\cD$ is the set of distinguished involutions, and $n_d=\pm1$ is as in Lusztig's property~$\mathsf{P5}$. 
\end{compactenum}
\smallskip

\noindent We note that the Losev-Bonnaf\'e left (respectively right) action corresponds to right (respectively left) multiplication by $\psi_0(\tw)$ in the above statements. 

Central to our strategy is a concrete construction of the elements~$\tw$ in Section~\ref{sec:concrete} that allows us to gain precise control over the denominators that are required for the element $\tw$ to exist. The next step is to study the action of $\tw$ on the Kazhdan-Lusztig basis in Section~\ref{sec:actiononKL}. In other words, we study the rational functions $b_{z,w}(\sq)\in\QQ(\sq)$ appearing in the expansion 
$$
\tw C_w=\sum_{z\in W}b_{z,w}(\sq)C_z. 
$$
It turns out that properties of the canonical trace $\Tr:\tilde{H}\to \QQ(\sq)$, along with the facts that $\tw^2=1$ and $\tw^*=\tw$, force considerable regularity of the coefficients $b_{z,w}(\sq)$. Moreover, using Bonnaf\'e's action of $T_{\sw_I}$ on the Kazhdan-Lusztig basis of $H$ from~\cite{Bon:16} we are able to determine one particular coefficient $b_{z,w}(\sq)$ precisely, and the regularity of the remaining coefficients, together with restrictions arising from the trace functional, force all other coefficients $b_{z,w}(\sq)$ to vanish at $\sq=0$. In Theorem~\ref{thm:mainresult} we state and prove our main theorem, and the Rouquier-White conjecture follows. We conclude the paper by explicitly computing the cactus group action in type $\tilde{\sA}_2$ in Section~\ref{sec:example}.

\begin{assumptions}\label{ass:1}
We will work in the setting of Kazhdan-Lusztig theory with general parameters (possibly unequal). In this setting there are various fundamental properties of cells, Lusztig's $\ba$-function, and the asymptotic algebra that are only conjectural. In \cite{Lus:03} Lusztig stated 15 properties $\mathsf{P1}$--$\mathsf{P15}$ that capture the desired (and expected) behaviour of the key ingredients of Kazhdan-Lusztig theory. In particular, these properties ensure the existence of the asymptotic algebra. Therefore, in the present paper, it is appropriate to assume properties $\mathsf{P1}$--$\mathsf{P15}$. We note that these properties are known to hold in equal parameters (the setting of the original Rouquier-White conjecture). Moreover, we will assume that the set $\cD$ of \textit{distinguished involutions} is finite -- this is implicitly assumed in~\cite{RW:24}, as they assume $\cJ$ is unital and that Lusztig's homomorphism $\psi:H\to \ZZ[\sq,\sq^{-1}]\otimes_{\ZZ}\cJ$ is defined. Thus our proof is unconditional in the important equal parameter case. 
\end{assumptions}

\section{Background and preliminaries}\label{sec:background}

In this section we recall background on Coxeter groups, cactus groups, the Hecke algebra, Kazhdan-Lusztig theory, the Mathas-Lusztig theorem, and Bonnaf\'e's action of the cactus group on~$W$.

\subsection{Coxeter groups and cactus groups}

Let $(W,S)$ be a Coxeter system. A subset $I\subseteq S$ is called \textit{spherical} if the group $W_I=\langle I\rangle$ is finite. If $I$ is spherical, let $\sw_I$ denote the longest element of $W_I$, and let $\omega_I:I\to I$ be given by $\omega_I(s)=\sw_Is\sw_I$. 

We say that a subset $I\subseteq S$ is \textit{connected} if the associated subgraph of the Coxeter graph of $(W,S)$ is connected (equivalently, $(W_I,I)$ is an irreducible Coxeter system).

\begin{defn}
The \textit{cactus group} $\cact(W)$ has generators $\gamma_I$ indexed by the connected spherical subsets $I\subseteq S$, with relations
\begin{align*}
\gamma_I^2&=1\\
\gamma_I\gamma_J&=\gamma_J\gamma_I&&\text{if $I\cap J=\emptyset$ and $st=ts$ for all $s\in I$ and $t\in J$}\\
\gamma_I\gamma_J\gamma_I&=\gamma_{\omega_I(J)}&&\text{if $J\subseteq I$}.
\end{align*}
\end{defn} 

\begin{remark}\label{rem:origins}
The origins of \textit{cactus groups} can be traced back to the work of Drinfeld~\cite{Dr:90}, Devadoss~\cite[Definition~6.1.2]{De:99} (as \textit{quasibraid groups}),  Davis, Januszkiewicz, Scott~\cite[Theorem~4.7.2]{DJS:03} (as \textit{mock reflection groups}), and Henriques and Kamnitzer~\cite{HK:06} (where the name \textit{cactus groups} is introduced). In the latter work, the action of the cactus group on tensor products in coboundary monoidal categories is analogous to the action of the braid group on tensor products in braided monoidal categories (see \cite[p.745]{EHKR:10}). Over the subsequent years a vast literature on cactus groups has developed, with connections to various branches of mathematics, see for example~\cite{BGL:19,CGP:20,EHKR:10,Gen:25,TK:09,LY:23}. 
\end{remark}

\subsection{The Hecke algebra}

Let $(W,S)$ be a Coxeter system, and let $\varphi:W\to \ZZ_{\geq 0}$ be a \textit{positive weight function}. That is, $\varphi(s)>0$ for all $s\in S$, and $\varphi(uv)=\varphi(u)+\varphi(v)$ whenever $\ell(uv)=\ell(u)+\ell(v)$. Let $H=H(W,S,\varphi)$ be the associated \textit{weighted Hecke algebra}. Thus $H$ is the algebra over $\sR=\ZZ[\sq,\sq^{-1}]$ with basis $(T_w\mid w\in W)$ and multiplication given by
$$
T_wT_s=\begin{cases}
T_{ws}&\text{if $\ell(ws)=\ell(w)+1$}\\
T_{ws}+(\sq^{\varphi(s)}-\sq^{-\varphi(s)})T_w&\text{if $\ell(ws)=\ell(w)-1$}.
\end{cases}
$$
In this definition, the parameter $\sq$ is an indeterminate. The Hecke algebra can be specialised at any nonzero real number. Specialising at $\sq=1$ the Hecke algebra becomes the group algebra $\ZZ[W]$ of~$W$. The ``equal parameter'' case is $\varphi=\ell$. 

\begin{remark}
The weight function $\varphi:W\to \ZZ_{\geq 0}$ will not play a significant role in our analysis and proofs, and thus it will be suppressed from the notation. However we emphasise that the structures we encounter later (the Kazhdan-Lusztig basis, the associated cells, and ultimately the action of the cactus group on the asymptotic algebra) all depend rather heavily on~$\varphi$. 
\end{remark}

If $I\subseteq S$, let $H_I$ be the subalgebra of $H$ spanned by $\{T_w\mid w\in W_I\}$. If $I$ is spherical then the automorphism $\omega_I:W_I\to W_I$ lifts to an algebra automorphism $\omega_I:H_I\to H_I$ with $\omega_I(T_w)=T_{\omega_I(w)}$. The following lemma is well known. 

\begin{lemma}\label{lem:omega}
If $I$ is spherical, then $T_{\sw_I}hT_{\sw_I}^{-1}=\omega_I(h)$ for all $h\in H_I$, and $\sT_{\sw_I}^2$ is central in~$H_I$. 
\end{lemma}

\begin{proof}
For the first claim, it suffices to check that $T_{\sw_I}T_sT_{\sw_I}^{-1}=T_{\omega_I(s)}$ for $s\in I$, and this is a basic consequence of the defining multiplication formulae in~$H_I$. The second statement follows as $\omega_I^2=1$ (since $\sw_I^2=1$). 
\end{proof}

The \textit{adjoint} anti-automorphism $*:H\to H$ is given by
$$
\big(\sum a_wT_w\big)^*=\sum a_wT_{w^{-1}}. 
$$
The \textit{canonical trace} on $H$ is the linear map $\Tr:H\to\ZZ[\sq,\sq^{-1}]$ given by (see \cite{GP:00})
$$
\Tr\big(\sum a_wT_w\big)=a_1. 
$$
We have $\Tr(T_u^*T_v)=\delta_{u,v}$, and it follows that $\Tr(h_1h_2)=\Tr(h_2h_1)$ for all $h_1,h_2\in H$. Thus we have a symmetric bilinear form given by
$$
\langle h_1,h_2\rangle=\Tr(h_1^*h_2),
$$
and if $h=\sum a_wT_w$ we write
$$
\|h\|^2=\langle h,h\rangle=\sum_{w\in W}a_w^2. 
$$
In particular, for all nonzero real specialisations of $\sq$ this bilinear form is positive definite. 

At times we will work in the algebra $\tilde{H}=\QQ(\sq)\otimes_{\sR} H$, where scalars are extended to the field of rational functions in~$\sq$. The adjoint map obviously extends to a map $*:\tilde{H}\to \tilde{H}$, and the canonical trace extends to a map $\Tr:\tilde{H}\to\QQ(\sq)$. 

\subsection{Kazhdan-Lusztig theory}\label{sec:KL}

We refer to \cite{Bon:17} and \cite{Lus:03} for general Kazhdan-Lusztig theory for Hecke algebras with unequal parameters. 

The \textit{bar-involution} on $H$ is the $\ZZ$-linear map with $\overline{\sq}=\sq^{-1}$ and $\overline{T_w}= T_{w^{-1}}^{-1}$. Let $(C_w\mid w\in W)$ be the \textit{Kazhdan-Lusztig basis}. Thus $C_w\in H$ is the unique element with
$$
\overline{C_w}=C_w\quad\text{and}\quad C_w=T_w+\sum_{v<w}P_{v,w}(\sq^{-1})T_v,
$$
where $<$ denotes Bruhat order, and $P_{v,w}(\sq^{-1})\in\sq^{-1}\ZZ[\sq^{-1}]$ for $v<w$ are the \textit{Kazhdan-Lusztig polynomials}. We write $P_{v,w}(\sq^{-1})=0$ if $v\not\leq w$, and $P_{w,w}(\sq^{-1})=1$. 

For $x,y\in W$ write 
\begin{align}\label{eq:h}
C_xC_y=\sum_{z\in W}h_{x,y,z}C_z
\end{align}
where $h_{x,y,z}\in\sR$ (necessarily $\overline{h_{x,y,z}}=h_{x,y,z}$). 

Let $\leq_L$, $\leq_R$, and $\leq_{LR}$ denote the left, right, and two-sided Kazhdan-Lusztig preorders, and let $\sim_L$, $\sim_R$, and $\sim_{LR}$ denote the associated equivalence relations on~$W$. Equivalence classes for these relations are left, right, and two-sided \textit{cells} (respectively). We have $x\leq_L y$ if and only if $x^{-1}\leq_R y^{-1}$, and hence $x\sim_L y$ if and only if $x^{-1}\sim_R y^{-1}$. 

We write $H^{\leq_R w}$ and $H^{<_R w}$ for the $\sR$-modules
$$
H^{\leq_R w}=\bigoplus_{v\leq_R w}\sR C_v\quad\text{and}\quad H^{<_Rw}=\bigoplus_{v<_R w}\sR C_v,
$$
where $v<_Rw$ means $v\leq_R w$ with $v\not\sim_R w$. We similarly define $\sR$-modules $H^{<_Lw}$, $H^{\leq_L w}$, $H^{<_{LR}w}$, and $H^{\leq_{LR}w}$. 

If $\Ga$ is a right cell then we write $H^{\leq_R \Ga}=H^{\leq_R w}$ and $H^{<_R\Ga}=H^{<_Rw}$ for any $w\in\Ga$, and the \textit{right cell module} of $\Ga$ is 
$$
M_{\Ga}=H^{\leq_R\Ga}/H^{<_R\Ga}. 
$$
By the definition of $\leq_R$, the cell module $M_{\Ga}$ is a right $H$-module with basis $(c_y^{\Ga}\mid y\in\Ga)$, where $c_y^{\Ga}=C_y+H^{<_R\Ga}$ is the image of $C_y$ in the quotient. There are similar constructions for left cell modules, and two-sided cell bimodules.

Let $\deg:\sR\to \ZZ\cup\{-\infty\}$ denote degree in~$\sq$. Since $\overline{h_{x,y,z}}=h_{x,y,z}$ we have $\deg h_{x,y,z}\geq 0$ whenever $h_{x,y,z}\neq 0$. \textit{Lusztig's $\ba$-function} is the function $\ba:W\to \ZZ_{\geq 0}$ defined by 
$$
\ba(z)=\max\{\deg h_{x,y,z}\mid x,y\in W\}. 
$$
The elements $\gamma_{x,y,z^{-1}}\in\ZZ$ are defined by 
$$
h_{x,y,z}=\gamma_{x,y,z^{-1}}\sq^{\ba(z)}+\text{(terms of degree $<\ba(z)$)}. 
$$
The set of \textit{distinguished involutions} is 
$$
\cD=\{z\in W\mid \ba(z)=-\deg P_{1,z}(\sq^{-1})\}. 
$$

Let $\mathsf{P1}$--$\mathsf{P15}$ be the properties listed in \cite[Conjecture~14.2]{Lus:03}. We list a subset of these properties below for convenience. 
\begin{compactenum}
\item[$\mathsf{P2}$.] If $d\in\cD$ and $x,y\in W$ satisfy $\gamma_{x,y,d}\neq 0$ then $y=x^{-1}$.
\item[$\mathsf{P4}$.] If $z'\leq _{LR}z$ then $\ba(z')\geq \ba(z)$. In particular, $\ba$ is constant on two-sided cells. 
\item[$\mathsf{P5}$.] If $d\in\cD$, $y\in W$, and $\gamma_{y,y^{-1},d}\neq 0$, then $\gamma_{y,y^{-1},d}=n_d\in\{-1,1\}$.
\item[$\mathsf{P6}$.] If $d\in\cD$ then $d^2=1$. 
\item[$\mathsf{P7}$.] For $x,y,z\in W$ we have $\gamma_{x,y,z}=\gamma_{y,z,x}$. 
\item[$\mathsf{P8}$.] If $\gamma_{x,y,z^{-1}}\neq 0$ then $z\sim_R x$, $z\sim_L y$, and $y\sim_R x^{-1}$. 
\item[$\mathsf{P9}$.] If $z'\leq_Lz$ and $\ba(z')=\ba(z)$ then $z'\sim_Lz$.
\item[$\mathsf{P12}$.] If $I\subseteq S$ then the $\ba$-function of $W_I$ is the restriction of the $\ba$-function of $W$.
\item[$\mathsf{P13}$.] Each right cell $\Ga$ contains a unique element $d\in\cD$, and $\gamma_{x,x^{-1},d}\neq 0$ for all $x\in\Ga$. 
\item[$\mathsf{P15}$.] If $u,v,z,z'\in W$ are such that $\ba(z)=\ba(z')$ then 
$$
\sum_{w\in W}h_{w,v,z}\otimes h_{u,z',w}=\sum_{w\in W}h_{z',v,w}\otimes h_{u,w,z}\quad\text{in $\sR\otimes_{\ZZ}\sR$}.
$$
\end{compactenum}

We will assume that properties $\mathsf{P1}$--$\mathsf{P15}$ hold throughout this paper (see Assumptions~\ref{ass:1}). Recall that $\mathsf{P1}$--$\mathsf{P15}$ are known to hold in the equal parameter case (see \cite[Section 15]{Lus:03}), which is the setting of the original Rouquier-White conjecture. Further instances where $\mathsf{P1}$--$\mathsf{P15}$ are known to hold include the ``quasisplit case'' \cite[Section 16]{Lus:03} (where the parameters arise from a diagram twist), finite type $\sF_4$ for all parameters~\cite{Geck:11}, affine types $\tilde{\sC}_2$ and $\tilde{\sG}_2$ for all parameters~\cite{GP:19,GP:19b}, and Coxeter groups with complete graph for all parameters~\cite{Xie:21}.

The \textit{asymptotic algebra} is the $\ZZ$-algebra $\cJ$ with basis $(\st_w\mid w\in W)$ and multiplication given by 
$$
\st_x\st_y=\sum_{z\in W}\gamma_{x,y,z^{-1}}\st_z. 
$$
Assuming $\mathsf{P1}$--$\mathsf{P15}$ it is known that $\cJ$ is an associative algebra~\cite[Section 18.3]{Lus:03}. If $\cD$ is finite then $\cJ$ is unital, with 
$$
1_{\cJ}=\sum_{d\in \cD}n_d\st_d,
$$
where $n_d=\gamma_{d,d,d}\in\{-1,1\}$. We note that $\cD$ is trivially finite for spherical Coxeter groups, and Lusztig proves in \cite[Theorem 18.2]{Lus:03} that (under the assumption $\mathsf{P1}$--$\mathsf{P15}$) $\cD$ is finite for affine Coxeter groups for all parameters. On the other hand, there are examples of universal Coxeter groups where $\cD$ is infinite (see \cite[Remark~4.5]{SY:15}). 

There is an anti-involution $*:\cJ\to\cJ$ given by $\st_w^*=\st_{w^{-1}}$. The \textit{canonical trace} on $\cJ$ (see \cite[Section 20.1]{Lus:03} and \cite[Section~5.1]{CGLP:24}) is the linear map $\Tr_{\infty}:\cJ\to\ZZ$ with 
$$
\Tr_{\infty}\big(\sum a_w\st_w\big)=\sum_{d\in\cD}n_da_d,
$$
and setting $\langle a,b\rangle_{\infty}=\Tr_{\infty}(a^*b)$ defines a positive definite symmetric bilinear form on~$\cJ$, with $\langle \st_x,\st_y\rangle_{\infty}=\delta_{x,y}$. Moreover, $\Tr_{\infty}(ab)=\Tr_{\infty}(ba)$ and $\langle ab,c\rangle_{\infty}=\langle b,a^*c\rangle_{\infty}$.

\subsection{Lusztig's homomorphism}\label{sec:LusHom}

Recall that we assume that $\cD$ is finite. In \cite[Theorem~18.9]{Lus:03} Lusztig defined a homomorphism of unital algebras $\phi:H\to\ZZ[\sq,\sq^{-1}]\otimes_{\ZZ}\cJ$ (assuming $\mathsf{P1}$--$\mathsf{P15}$). For our purposes (and following the setup of \cite{RW:24} and \cite[Chapter~19]{Bon:17}), it is convenient to consider the algebra homomorphism $\psi=\phi\circ\dagger$, where $\dagger:H\to H$ is the involution $T_w^{\dagger}=(-1)^{\ell(w)}T_{w^{-1}}^{-1}$. Thus we have an algebra homomorphism $
\psi:H\to \ZZ[\sq,\sq^{-1}]\otimes_{\ZZ}\cJ$ with
$$
\psi(C_x)=\sum_{d\in\cD,\,y\in W,\,\ba(y)=\ba(d)}n_dh_{x,d,y}\st_y.
$$
Note that the definition of $\psi$ requires $\cD$ to be finite (otherwise $\psi(C_1)$ is an infinite sum). Moreover, note that the sum over $y\in W$ does not introduce infinitely many nonzero terms, because for each fixed $x\in W$ and $d\in\cD$ there are only finitely many $y\in W$ with $h_{x,d,y}\neq 0$ (from the definition~(\ref{eq:h})).

A main property of $\psi$ that we will require is the following observation (cf. \cite[Theorem~5.2]{Geck:09}). 

\begin{lemma}\label{lem:51}
Let $h\in H$ and $w\in W$. If
$$
hC_w=\sum_z b_{z,w}C_z\quad\text{then}\quad 
\psi(h)\st_w=\sum_{z\sim_{LR}w}b_{z,w}\st_z.
$$
\end{lemma}

\begin{proof}
By linearity, it suffices to consider $h=C_x$, where we have $C_xC_w=\sum_zh_{x,w,z}C_z$. On the other hand, we have
\begin{align*}
\psi(C_x)\st_w&=\sum_{d\in\cD,\,y\in W,\,\ba(y)=\ba(d)}n_dh_{x,d,y}\st_y\st_w\\
&=\sum_z\bigg(\sum_{d\in\cD,\,y\in W,\,\ba(y)=\ba(d)}n_dh_{x,d,y}\gamma_{y,w,z^{-1}}\bigg)\st_z. 
\end{align*}
A nonzero summand satisfies $\gamma_{y,w,z^{-1}}\neq 0$ and hence by $\mathsf{P8}$ we have $y\sim_L w^{-1}$, $w\sim_L z$, and $z^{-1}\sim_L y^{-1}$. Thus by $\mathsf{P4}$ a nonzero summand satisfies $\ba(d)=\ba(y)=\ba(w)=\ba(z)$. Now $\mathsf{P15}$ (with $(u,v,w,z,z')=(x,w,y,z,d)$) gives the identity 
$$
\sum_y \gamma_{y,w,z^{-1}}h_{x,d,y}=\sum_y\gamma_{d,w,y^{-1}}h_{x,y,z}
$$
when $\ba(d)=\ba(w)=\ba(z)$ (see the proof of \cite[Theorem~18.9]{Lus:03}), and so in the formula for $\psi(C_x)\st_w$, the coefficient of $\st_z$, for $\ba(z)=\ba(w)$, becomes
$$
\sum_y\bigg(\sum_{d\in\cD}n_d\gamma_{d,w,y^{-1}}\bigg)h_{x,y,z}. 
$$
Using $\mathsf{P7}$ and $\mathsf{P2}$ we have $\sum_y\left(\sum_{d\in\cD}n_d\gamma_{d,w,y^{-1}}\right)=\sum_y\left(\sum_{d\in\cD}n_d\gamma_{w,y^{-1},d}\right)=\sum_{d\in\cD}n_d\gamma_{w,w^{-1},d}$. By $\mathsf{P13}$ and $\mathsf{P5}$ we have $\gamma_{w,w^{-1},d}=0$ if $d\not\sim_R w$ and $\gamma_{w,w^{-1},d}=n_d$ if $d\sim_R w$, and so the coefficient of $\st_z$ reduces to just $h_{x,w,z}$. So 
$$
\psi(C_x)\st_w=\sum_{\ba(z)=\ba(w)}h_{x,w,z}\st_z. 
$$
Finally, if $h_{x,w,z}\neq 0$ then $z\leq_L w$. By $\mathsf{P4}$ this gives $\ba(z)\geq \ba(w)$, and since equality holds $\mathsf{P9}$ gives $z\sim_L w$, hence $z\sim_{LR}w$. Thus the sum above is over $z\sim_{LR}w$. 
\end{proof}

\subsection{The Mathas-Lusztig Theorem}\label{sec:ML}

If $I$ is spherical we write (for $w\in W_I$)
$$
\alpha_I(w)=\ba(\sw_Iw)-\ba(w).
$$
Note that $\alpha_I$ is constant on each right cell $\Ga$ of $W_I$, for if $x,y\in W_I$ with $x\leq_R y$ then $\sw_I y\leq_R\sw_I x$ (see \cite[Corollary~11.7]{Lus:03}). In particular, if $x\sim_R y$ then $\sw_I x\sim_R\sw_I y$. Then by $\mathsf{P4}$ it follows that for $x\sim_R y$ we have $\alpha_I(x)=\ba(\sw_Ix)-\ba(x)=\ba(\sw_Iy)-\ba(y)=\alpha_I(y)$. Thus we can denote the constant value of $\alpha_I$ on $\Ga$ by $\alpha_I(\Ga)$. We note that if $\ba_I:W_I\to\ZZ_{\geq 0}$ is the $\ba$-function of $H_I$, then $\ba_I$ is the restriction of $\ba:W\to \ZZ_{\geq 0}$ to $W_I$, by~$\mathsf{P12}$.

The following theorem was first proved by Mathas~\cite{Mat:96} in the equal parameter case, and extended by Lusztig~\cite[Theorem~2.3]{Lus:15} in unequal parameters. 

\begin{thm}[Mathas-Lusztig]\label{thm:ML}
Let $I\subseteq S$ be spherical. There exists a unique sign map $\eta_I:W_I\to \{-1,1\}$, and unique involutions $\rho_I:W_I\to W_I$ and $\lambda_I:W_I\to W_I$ such that, for all $w\in W_I$, 
\begin{align*}
\sq^{\alpha_I(w)}T_{\sw_I}C_w&=\eta_I(w)C_{\rho_I(w)}\pmod{H_I^{<_{LR}w}}\\
\sq^{\alpha_I(w)}C_wT_{\sw_I}&=\eta_I(w)C_{\lambda_I(w)}\pmod{H_I^{<_{LR}w}}.
\end{align*}
We have $\la_I(w)=\rho_I(w^{-1})^{-1}$, $\rho_I=\la_I\circ\omega_I$, and $\rho_I(w)\sim_L w$ and $\la_I(w)\sim_R w$ for all $w\in W_I$. 
\end{thm}

\begin{remark}\label{rem:etaequal}
In the equal parameter case the sign map $\eta_I$ is given by $\eta_I(w)=(-1)^{\ba(\sw_Iw)}$. See \cite[Example~17.3.2]{Bon:16} (noting that $\eta_I(w)=\eta_I(w^{-1})$). 
\end{remark}

\subsection{Parabolic projections and parabolic lifts}

Let $I\subseteq S$. Each right (respectively left) $W_I$-coset has a unique minimal length element. Let ${^I}W$ (respectively $W^I$) denote the set of these minimal length elements. Each element $w\in W$ can be written in a unique way as $w=xy=y'x'$ with $x,x'\in W_I$, $y\in {^I}W$, and $y'\in W^I$. The \textit{right parabolic projection} of $w$ onto $W_I$ is 
$$
\pr_I^R(w)=x,
$$
and similarly we define the \textit{left parabolic projection} of $w$ onto $W_I$ by
$$
\pr_I^L(w)=x'. 
$$
Following Bonnaf\'e~\cite{Bon:16}, if $f:W_I\to X$ is a function (with $X$ any set), we define functions $f_R:W\to X$ and $f_L:W\to X$ by 
$$
f_R(w)=f(\pr_I^R(w))\quad\text{and}\quad f_L(w)=f(\pr_I^L(w)). 
$$
If $f:W_I\to W_I$ is a function, we define functions $f^R:W\to W$ and $f^L:W\to W$ by 
$$
f^R(w)=f(x)y\quad\text{and}\quad f^L(w)=y'f(x'),
$$
where $w=xy=y'x'$ with $x,x'\in W_I$, $y\in {^I}W$, and $y'\in W^I$ (thus $x=\pr_I^R(w)$, $x'=\pr_I^L(w)$). Finally, if $f:W\to W$ we write $f^{\mathsf{op}}:W\to W$ for the function $f^{\mathsf{op}}(w)=f(w^{-1})^{-1}$.

The above notation will be used, in particular, for the sign functions $\eta_I:W_I\to \{-1,1\}$, the function $\alpha_I:W_I\to\ZZ$, and the involutions $\la_I:W_I\to W_I$ and $\rho_I:W_I\to W_I$ from the Mathas-Lusztig Theorem. Thus if $w=xy$ with $x\in W_I$ and $y\in {^I}W$ (so $\pr_I^R(w)=x$) then 
\begin{align*}
\alpha_{I,R}(w)=\alpha_I(x),\quad \eta_{I,R}(w)=\eta_I(x),\quad\text{and}\quad \rho_I^R(w)=\rho_I(x)y
\end{align*}
and similarly for $\alpha_{I,L}(w)$, $\eta_{I,L}(w)$, and $\lambda_I^L(w)$ (using left projections).

\subsection{Bonnaf\'{e}'s induced action}

In \cite[Corollary 4.10]{Bon:16} Bonnaf\'e lifts the action of $T_{\sw_I}$ on the spherical Hecke algebra $H_I$ (from the Mathas-Lusztig Theorem) to the full ambient Hecke algebra. In particular, he proves the following. 

\begin{thm}[{\cite[Corollaries~4.7 and 4.10]{Bon:16}}]\label{thm:bonnafe}
Let $I\subseteq S$ be spherical. For $w\in W$ we have
\begin{align*}
\sq^{\alpha_{I,R}(w)}T_{\sw_I}C_w=\eta_{I,R}(w)C_{\rho_I^R(w)}\pmod{H_I^{<_R\omega_{I,R}(w)}H}\\
\sq^{\alpha_{I,L}(w)}C_wT_{\sw_I}=\eta_{I,L}(w)C_{\la_I^L(w)}\pmod{HH_I^{<_L\omega_{I,L}(w)}}
\end{align*}
Moreover, 
\begin{compactenum}[$(1)$]
\item $w\sim_R\la_I^L(w)$ and $w\sim_L\rho_I^R(w)$ for all $w\in W$;
\item if $x\sim_L y$ then $\la_I^L(x)\sim_L\la_I^L(y)$;
\item if $x\sim_R y$ then $\rho_I^R(x)\sim_R \rho_I^R(y)$;
\item $(\la_I^L)^{\mathsf{op}}=\rho_I^R$, and $\eta_{I,L}(w)=\eta_{I,R}(w^{-1})$. 
\end{compactenum}
\end{thm}

\begin{remark}\label{rem:bonnafeaction}
Bonnaf\'e's action of $\cact(W)\times\cact(W)$ on $W$ given in the introduction is given by $\gamma_I^L=\la_I^L$ and $\gamma_I^R=\rho_I^R$. 
\end{remark}

\section{Proof of main theorem}

In this section we give the proof of~\cite[Conjecture~4.2]{RW:24}. As outlined in the introduction, Section~\ref{sec:concrete} gives our construction of the elements $\tw\in\QQ(\sq)\otimes_{\sR}H$, Section~\ref{sec:actiononKL} determines information about the action of $\tw$ on the Kazhdan-Lusztig basis, and the statement (and proof) of our main theorem is given in Section~\ref{sec:mainproof}. 

\subsection{Construction of the element $\tw$}\label{sec:concrete}

The elements $\tw$, for $I\subseteq S$ spherical, are constructed in \cite{EHKR:10,RW:24} using ``Drinfeld's unitarisation trick''~\cite{Dr:90}, and pass through a completion of the braid group. It is convenient for our purposes to give a concrete construction of $\tw$ directly in the Hecke algebra. While this construction is less efficient, it allows us to control more precisely the extension of the coefficient ring that is required for the existence of the elements (see in particular Corollary~\ref{cor:coefficientring}).

We first give some motivation for the elements $\tw$. We will show in Theorem~\ref{thm:unique} below that for each spherical $I\subseteq S$ there exists a unique element $\tw\in\QQ(\sq)\otimes H_I$ such that $\tw^2=1$, $\tw h\tw^{-1}=\omega_I(h)$ for all $h\in\QQ(\sq)\otimes_{\sR} H_I$, and $\tw$ specialises at $\sq=1$ to $\sw_I$. Assuming this for the moment, we have the following elementary observation.

\begin{prop}\label{prop:action} The assignment $\gamma_I\cdot h=\tw h$ for $h\in \QQ(\sq)\otimes_{\sR}H$ defines a linear action of the cactus group $\cact(W)$ on the Hecke algebra $\QQ(\sq)\otimes_{\sR}H$. In other words, the map $\cact(W)\to (\QQ(\sq)\otimes_{\sR}H)^{\times}$ with $\gamma_I\mapsto \tw$ is a group homomorphism. 
\end{prop}

\begin{proof} 
By Theorem~\ref{thm:unique} we have $\tw^2=1$. If $I,J\subseteq S$ are connected spherical subsets with $J\subseteq I$ then property (2) in Theorem~\ref{thm:unique} (along with uniqueness) gives
$
\tw\widetilde{\sw}_J\tw=\widetilde{\sw}_{\omega_I(J)}. 
$
On the other hand, if $I\cap J=\emptyset$ and $st=ts$ for all $s\in I$ and $t\in J$ then $h_1h_2=h_2h_1$ for all $h_1\in \QQ(\sq)\otimes_{\sR}H_I$ and $h_2\in \QQ(\sq)\otimes_{\sR}H_J$, and hence $\tw\widetilde{\sw}_J=\widetilde{\sw}_J\tw$. 
\end{proof}

We now turn our attention to constructing the elements $\tw$ (and proving Theorem~\ref{thm:unique}). Recall from Section~\ref{sec:ML} that $\alpha_I(\Ga)$ denotes the constant value of $\alpha_I(w)=\ba(\sw_Iw)-\ba(w)$ on the right cell~$\Ga$. 

\begin{lemma}\label{lem:eigenvalues}
The characteristic polynomial of $T_{\sw_I}^2$, considered as an operator acting on $H_I$ on the right, is 
$$
\det(X -T_{\sw_I}^2)=\prod_{\Ga}(X-\sq^{-2\alpha_I(\Ga)})^{|\Ga|},
$$
where the product is over all right cells in $W_I$. In particular, every eigenvalue of $T_{\sw_I}^2$ is of the form $\sq^{2m}$ for some $m\in\ZZ$. Moreover $T_{\sw_I}^2$ is diagonalisable over~$\QQ(\sq)$.
\end{lemma}

\begin{proof}
Consider a right cell $\Ga$ of $W_I$. Let $M_{\Ga}$ be its right cell module, with Kazhdan-Lusztig basis $(c_y^{\Ga}\mid y\in\Ga)$. The Mathas-Lusztig Theorem (Theorem~\ref{thm:ML}) gives
$$
\sq^{\alpha_I(\Ga)}c_y^{\Ga}T_{\sw_I}=\eta_I(y)c_{\la_I(y)}^{\Ga}\quad\text{for $y\in \Ga$}.
$$
Thus right multiplication by $\sq^{\alpha_I(\Ga)}T_{\sw_I}$ on the module $M_{\Ga}$ is represented by a signed permutation matrix~$P_{\Ga}$, whose entries are independent of $\sq$. At $\sq=1$ the same matrix represents multiplication by $\sw_I$ on the specialised module, and since $\sw_I^2=1$ we have $P_{\Ga}^2=1$. Thus the right action of $T_{\sw_I}^2$ on $M_{\Ga}$ is given by the scalar $\sq^{-2\alpha_I(\Ga)}$. Since $H_I$ has a filtration of its right regular module with successive quotients given by right cell-modules, by ordering the Kazhdan-Lusztig basis elements $(C_y\mid y\in W_I)$ in a way compatible with the partial order on two-sided cells it follows that the matrix representing $T_{\sw_I}^2$ is upper triangular, and that the entry on the diagonal corresponding to $C_y$ is $\sq^{-2\alpha_I(\Ga)}$ if $y\in\Ga$. Hence the characteristic polynomial of $T_{\sw_I}^2$ is as claimed. 

Let $k$ be a characteristic~$0$ splitting field for $W_I$. Then $k(\sq)\otimes_{\sR}H_I$ is split semisimple (by the Tits Deformation Theorem~\cite[Theorem~7.4.6]{GP:00}). Since $T_{\sw_I}^2$ is central, it is diagonalisable over $k(\sq)$, however since all eigenvalues already lie in the field $\QQ(\sq)$ we conclude that $T_{\sw_I}^2$ is diagonalisable over~$\QQ(\sq)$. 
\end{proof}

Recall that if $T$ is a diagonalisable operator acting on a finite dimensional vector space over a field, and if $\la_1,\ldots,\la_r$ are the distinct eigenvalues of $T$, then 
$$
T=\sum_{i=1}^r\lambda_iE_i\quad\text{where}\quad E_i=\prod_{j\neq i}\frac{T-\la_jI}{\la_i-\la_j}
$$
Here $E_i$ is the projection onto the $\lambda_i$-eigenspace, and we have $E_iE_j=\delta_{i,j}E_i$ and $\sum E_i=1$. Moreover, if square roots $\la_i^{1/2}$ are chosen in the field (or a relevant field extension) then $T^{1/2}=\sum_{i=1}^r\la_i^{1/2}E_i$ satisfies $(T^{1/2})^2=T$.

\begin{defn}\label{defn:wI}
Let $\sq^{2m_1},\ldots,\sq^{2m_r}$ be the distinct eigenvalues of $T_{\sw_I}^2$ (see Lemma~\ref{lem:eigenvalues}), and for $1\leq i\leq r$ let 
$$
E_i=\prod_{j\neq i}\frac{T_{\sw_I}^2-\sq^{2m_j}}{\sq^{2m_i}-\sq^{2m_j}}\in \QQ(\sq)\otimes_{\sR}H_I
$$
be the projection onto the $\sq^{2m_i}$-eigenspace. Define $\tw\in\QQ(\sq)\otimes_{\sR} H_I$ by 
$$
\tw=T_{\sw_I}Z_I\quad\text{where}\quad Z_I=\sum_{i=1}^r\sq^{-m_i}E_i.
$$
\end{defn}

\begin{lemma}\label{lem:basic}
We have the following:
\begin{compactenum}[$(1)$]
\item $Z_I$ is central in $\QQ(\sq)\otimes_{\sR}H_I$, and $Z_I^2=T_{\sw_I}^{-2}$;
\item if $T_{\sw_I}^2$ acts on a right cell module (in $H_I$) by $\sq^{2m}$, then $Z_I$ acts on the module by $\sq^{-m}$;
\item $\tw^2=1$;
\item $\tw h\tw^{-1}=\omega_I(h)$ for all $h\in \QQ(\sq)\otimes_{\sR}H_I$;
\item $\tw$ is bar-invariant;
\item $\tw^*=\tw$. 
\end{compactenum}
\end{lemma}

\begin{proof}
(1) Since $Z_I$ is a polynomial in the central element $T_{\sw_I}^2$ it is itself central. We have $T_{\sw_I}^2Z_I^2=\left(\sum \sq^{2m_i}E_i\right)\left(\sum \sq^{-m_i}E_i\right)^2=\sum E_i=1$. 

(2) If $T_{\sw_I}^2$ acts by $\sq^{2m}$ then only the spectral projector $E_i$ with $m_i=m$ acts nontrivially on the module, and so $Z_I$ acts by $\sq^{-m_i}$. 

(3) We have $\tw^2=Z_I^2T_{\sw_I}^2=1$, by (1). 

(4) Since $Z_I$ is central, we have $\tw h\tw^{-1}=T_{\sw_I}hT_{\sw_I}^{-1}=\omega_I(h)$, by Lemma~\ref{lem:omega}. 

(5) Applying the bar involution to $E_iT_{\sw_I}^2=\sq^{2m_i}E_i$ gives $\overline{E_i}T_{\sw_I}^{-2}=\sq^{-2m_i}\overline{E_i}$, and by the uniqueness of the spectral projector it follows that $\overline{E_i}=E_i$. Therefore 
$$\overline{Z_I}=\sum_{i=1}^r \sq^{m_i}E_i=\sum_{i=1}^r \sq^{2m_i}\sq^{-m_i}E_i=T_{\sw_I}^2Z_I.
$$ 
It follows that 
$
\overline{\tw}=\overline{T_{\sw_I}}\,\overline{Z_I}=T_{\sw_{I}}^{-1}T_{\sw_I}^2Z_I=T_{\sw_I}Z_I=\tw,
$
hence the result. 

(6) By the formula in Definition~\ref{defn:wI} we have $E_i^*=E_i$ (as $T_{\sw_I}^*=T_{\sw_I}$), and so $\tw^*=\tw$. 
\end{proof}

We say that a rational function $f(\sq)\in\QQ(\sq)$ is \textit{regular} at $\alpha\in\RR$ if it has no pole at $\sq=\alpha$. Equivalently, if $f(\sq)=P(\sq)/Q(\sq)$ in reduced form, then $Q(\alpha)\neq 0$.

\begin{lemma}\label{lem:regular}
For $1\leq i\leq r$ let 
$$
E_i=\sum_{x\in W_I}a_{i,x}(\sq)T_x,\quad\text{with}\quad a_{i,x}(\sq)\in\QQ(\sq). 
$$
For each $x\in W_I$ the rational function $a_{i,x}(\sq)$ is regular for all real values of $\sq$. In particular, the element $E_i$ admits a specialisation for all nonzero real numbers~$\sq$. 
\end{lemma}

\begin{proof}
Fix $1\leq i\leq r$, and let $q\neq 0$ be a real number such that $E_i$ admits a specialisation $E_i(q)$ at $\sq=q$. Note that $E_i^*=E_i$ (as $E_i$ is a polynomial in the self-adjoint element $T_{\sw_I}^2$). Thus $1-E_i$ is a self-adjoint idempotent element and so on specialising
$$
\Tr(1-E_i(q))=\Tr((1-E_i(q))^*(1-E_i(q)))=\|1-E_i(q)\|^2\geq 0,
$$
and thus $\Tr(E_i(q))\leq 1$. Since
$
\Tr(E_i)=\Tr(E_i^*E_i)=\sum_{x\in W_I}a_{i,x}(\sq)^2
$, after specialising we have 
$$
\sum_{x\in W_I}a_{i,x}(q)^2\leq 1. 
$$
If there were a $q_0\in\RR$ such that $a_{i,x}(\sq)$ is not regular at $\sq=q_0$, then the above inequality shows that $a_{i,x}(q)$ is bounded for $q$ in a punctured real neighbourhood of $q_0$, a contradiction. Hence the result. 
\end{proof}

\begin{cor}\label{cor:specialise}
We have 
$$
\tw=\sum_{x\in W_I}a_x(\sq)T_x
$$
where each $a_x(\sq)$, $x\in W_I$, is regular at all real values of $\sq$. In particular, $\tw$ admits a specialisation for all nonzero real numbers. Moreover, $\tw$ specialises at $\sq=1$ to the element $\sw_I$ of the group algebra. 
\end{cor}

\begin{proof}
Since $\tw^*=\tw$ and $\tw^2=1$ (see Lemma~\ref{lem:basic}) we have 
$$
1=\Tr(\tw^*\tw)=\sum_{x\in W_I}a_x(\sq)^2,
$$
which shows that each $a_x(\sq)$ is regular at all real values of $\sq$. (Note that regularity at all nonzero values follows already from the definition of $\tw$ and Lemma~\ref{lem:regular}, however the above argument is needed to check the case $\sq=0$ as $\sq^{-1}$ may occur as eigenvalues in the definition of $\tw$, and multiplication by $T_{\sw_I}$ introduces $\sq^{-1}$ through the defining relations of $H_I$). 

It follows that $\tw$ admits a specialisation for all nonzero real numbers. Moreover, since $\sum_i E_i=1$, the element $Z_I$ specialises to $1$ at $\sq=1$, and thus $\tw=T_{\sw_I}Z_I$ specialises to $\sw_I$ at~$\sq=1$.
\end{proof} 

\begin{thm}\label{thm:unique}
For each spherical subset $I\subseteq S$ there exists a unique element $\tw\in\QQ(\sq)\otimes_{\sR}H_I$ such that 
\begin{compactenum}[$(1)$]
\item $\tw^2=1$;
\item $\tw h\tw^{-1}=\omega_I(h)$ for all $h\in\QQ(\sq)\otimes_{\sR}H_I$;
\item $\tw$ specialises at $\sq=1$ to the element $\sw_I$. 
\end{compactenum}
\end{thm}

\begin{proof}
Existence is given by Lemma~\ref{lem:basic} and Corollary~\ref{cor:specialise}. To prove uniqueness we argue as in \cite[Proposition~3.2]{RW:24}. Suppose that $a,b\in\QQ(\sq)\otimes_{\sR}H_I$ both satisfy (1), (2) and (3), and let $z=ab$. Then by (2) and the fact that $\omega_I^2=1$ we have that $z$ is central in $\QQ(\sq)\otimes_{\sR}H_I$. Thus $z^2=abz=azb=a^2b^2=1$, by (1). So $(z-1)(z+1)=0$. Working over $k(\sq)\otimes_{\sR}H_I$ (with $k$ a characteristic~$0$ splitting field for $W_I$ as in the proof of Lemma~\ref{lem:eigenvalues}) the central element $z$ is diagonalisable. By (3) $z$ specialises to $1$ at $\sq=1$, and hence $z+1$ is invertible in $k(\sq)\otimes_{\sR}H_I$. Thus the equation $(z-1)(z+1)=0$ gives $z=1$, proving uniqueness. 
\end{proof}

The uniqueness in Theorem~\ref{thm:unique} implies that the elements $\tw$ constructed above agree with the elements, also denoted $\tw$, that are constructed in \cite[Proposition~3.2]{RW:24} in equal parameters.

\begin{example}
If $I=\{s\}$ then $\tw=\widetilde{\sw}_s$ is given by
$$
\widetilde{w}_s=\frac{1-\sq^{2a}}{1+\sq^{2a}}+\frac{2\sq^{a}}{1+\sq^{2a}}T_s,
$$
where $a=\varphi(s)$ (one can directly check the defining properties; see also \cite[Proposition~3.4]{RW:24}). For another example, if $I=\{s,t\}$ and $W_I$ is a dihedral group of order~$8$ with $\varphi(s)=a>0$ and $\varphi(t)=b>0$ then writing $D=(\sq^{2a}+\sq^{2b})(1+\sq^{2a+2b})$ we have
$$
\tw=\frac{D-4\sq^{2a+2b}}{D}T_1+\frac{2\sq^{a+2b}(1-\sq^{2a})}{D}(T_s+T_{tst})+\frac{2\sq^{2a+b}(1-\sq^{2b})}{D}(T_t+T_{sts})+\frac{4\sq^{2a+2b}}{D}T_{stst}.
$$
A direct calculation in $\mathsf{MAGMA}$ verifies that this element satisfies the defining properties from Theorem~\ref{thm:unique}. 
\end{example}

The above examples illustrate that one must extend scalars to define the elements~$\tw$. Our construction gives the following precise control over the denominators that are required. 

\begin{cor}\label{cor:coefficientring}
With $a_x(\sq)$ as in Corollary~\ref{cor:specialise}, we have 
$$
a_x(\sq)\in\sR_I:=\ZZ[\sq,(1+\sq^2+\cdots+\sq^{2m-2})^{-1}\mid m\in M_I],
$$
where $M_I$ is a finite set of positive integers. Thus we have 
$$
\tw\in \bigoplus_{w\in W_I}\sR_IT_w\subseteq \tilde{H}. 
$$ 
\end{cor}

\begin{proof}
Since
$\sq^{2m_i}-\sq^{2m_j}=u(1-\sq^{2|m_i-m_j|})$ with $u$ a unit in $\ZZ[\sq,\sq^{-1}]$, the formula in Definition~\ref{defn:wI} already shows that for each $x\in W_I$,
$$
a_x(\sq)\in\ZZ[\sq,\sq^{-1},(1-\sq^{2m})^{-1}\mid m\in M_I],
$$
where $M_I$ is a finite set of positive integers. Write $1-\sq^{2m}=(1-\sq^2)(1+\sq^2+\cdots+\sq^{2m-2})$ (note that the latter factor has no real roots). Corollary~\ref{cor:specialise} implies that the denominator of the rational function $a_x(\sq)$ (when written in reduced form) cannot contain the factors $\sq$, $1-\sq$, or $1+\sq$, and since division by monic polynomials preserves integrality, it follows that 
$
a_x(\sq)\in\sR_I
$.
\end{proof}

\subsection{The action of $\tw$ on the Kazhdan-Lusztig basis}\label{sec:actiononKL}

We now turn attention to the action of $\tw$ on the Kazhdan-Lusztig basis of the ambient Coxeter group~$W$. The theorem below gives an extension of Bonnaf\'e's action from Theorem~\ref{thm:bonnafe}. In this section we work in the Hecke algebra over $\QQ(\sq)$. To abbreviate notation we write $\tilde{H}=\QQ(\sq)\otimes_{\sR}H$, and similarly for $\tilde{H}_I$. 

\begin{thm}\label{thm:actiononKLbasis}
Let $I\subseteq S$ be spherical. For $w\in W$ we have
\begin{align*}
\tw C_w=\eta_{I,R}(w)C_{\rho_I^R(w)}\pmod{\tilde{H}_I^{<_R\omega_{I,R}(w)}\tilde{H}}\\
C_w\tw=\eta_{I,L}(w)C_{\la_I^L(w)}\pmod{\tilde{H}\tilde{H}_I^{<_L\omega_{I,L}(w)}}
\end{align*}
\end{thm}

\begin{proof}
It is sufficient to prove the first equation (the second equation follows using the adjoint anti-involution). Let $w=xy$ with $x\in W_I$ and $y\in {^I}W$, and so $x=\pr_I^R(w)$. We must prove that
\begin{align}\label{eq:toprove}
\tw C_w=\eta_I(x)C_{\rho_I(x)y}\pmod{\tilde{H}_I^{<_R\omega_I(x)}\tilde{H}}.
\end{align}
Let $\Ga$ be the right cell of $W_I$ containing~$x$. From the proof of Lemma~\ref{lem:eigenvalues}, the element $T_{\sw_I}^2$ acts on the right-cell module $M_{\Ga}=\tilde{H}_I^{\leq_R x}/\tilde{H}_I^{<_R x}$ as the scalar $\sq^{-2\alpha_I(\Ga)}$, and so by Lemma~\ref{lem:basic} the element $Z_I$ acts on $M_{\Ga}$ as the scalar $\sq^{\alpha_I(\Ga)}$. Since $Z_I$ is central in $\tilde{H}_I$ we have 
$$
(Z_I-\sq^{\alpha_I(\Ga)})\tilde{H}_I^{\leq_Rx}\subseteq \tilde{H}_I^{<_Rx}
$$
under left multiplication as well, and hence 
\begin{align}\label{eq:diff}
(Z_I-\sq^{\alpha_I(\Ga)})\tilde{H}_I^{\leq_Rx}\tilde{H}\subseteq \tilde{H}_I^{<_Rx}\tilde{H}. 
\end{align}
From the right handed version of Geck's Induction Theorem (\cite{Ge:03}, see also \cite[Theorem~1.3]{Bon:16}) we have 
$$
\tilde{H}_I^{\leq_Rx}\tilde{H}=\bigoplus_{\substack{x'\in W_I\text{ with }x'\leq_Rx\\
y'\in {^I}W}}\QQ(\sq)C_{x'y'},
$$
and so in particular $C_w\in \tilde{H}_I^{\leq_Rx}\tilde{H}$. Hence by~(\ref{eq:diff}) we have
\begin{align}\label{eq:part1}
Z_IC_w=\sq^{\alpha_I(\Ga)}C_w\pmod{\tilde{H}_I^{<_Rx}\tilde{H}}. 
\end{align}

Now, since $T_{\sw_I}h=\omega_I(h)T_{\sw_I}$ for $h\in\tilde{H}_I$, and since $\omega_I$ preserves the right preorder, we have 
$$
T_{\sw_I}\tilde{H}_I^{<_R x}\tilde{H}=\tilde{H}_I^{<_R\omega_I(x)}\tilde{H},
$$
and~(\ref{eq:toprove}) now follows from Bonnaf\'e's Theorem (Theorem~\ref{thm:bonnafe}) and~(\ref{eq:part1}).
\end{proof}

\begin{remark}\label{rem:nokill}
Note that $C_{\rho_I^R(w)}\notin \tilde{H}_I^{<_R\omega_{I,R}(w)}\tilde{H}$, and hence the term $\eta_{I,R}(w)C_{\rho_I^R(w)}$ in Theorem~\ref{thm:actiononKLbasis} cannot be cancelled by an error term (and similarly for the left handed version). To see this, recall that $\rho_I=\la_I\circ\omega_I$. Writing $w=xy$ with $x\in W_I$ and $y\in {^I}W$ we have $\rho_I(x)\sim_R\la_I(\omega_I(x))\sim_R\omega_I(x)$ (since $\la_I(u)\sim_R u$ for all $u\in W_I$). Moreover, $\pr_I^R(\rho_I(x)y)=\rho_I(x)$. The module $\tilde{H}_I^{<_R\omega_I(x)}\tilde{H}$ is spanned by those $C_z$ for which $\pr_I^R(z)<_R\omega_I(x)$, and so in particular it does not contain the basis term $C_{\rho_I(x)y}$. 
\end{remark}

For $w\in W$, and $I\subseteq S$ spherical, consider the expansion
\begin{align}\label{eq:expand}
\tw C_w=\sum_{z\in W}b_{z,w}(\sq)C_z\quad\text{where $b_{z,w}(\sq)\in\QQ(\sq)$.}
\end{align}

\begin{thm}\label{thm:specialisation}
The coefficients $b_{z,w}(\sq)$ are regular at all real values, and more precisely
$$
b_{z,w}(\sq)\in\ZZ[\sq,(1+\sq^2+\cdots+\sq^{2m-2})^{-1}\mid m\in M_I],
$$
where $M_I$ is the finite set of integers in Corollary~\ref{cor:coefficientring}. Moreover, we have
$$
b_{z,w}(0)=\begin{cases}
\eta_{I,R}(w)&\text{if $z=\rho_I^R(w)$}\\
0&\text{otherwise}.
\end{cases}
$$
\end{thm}

\begin{proof}
It follows from Corollary~\ref{cor:coefficientring}, the definition of $C_w$, and the multiplication rules in the Hecke algebra that 
$$
b_{z,w}(\sq)\in\ZZ[\sq,\sq^{-1},(1+\sq^2+\cdots+\sq^{2m-2})^{-1}\mid m\in M_I]
$$
(note that we include $\sq^{-1}$ here). In particular, each rational function $b_{z,w}(\sq)$ is regular at all nonzero real values of $\sq$.

We claim that each rational function $b_{z,w}(\sq)$ is bounded as $\sq\to\infty$. Equivalently, we must show that $\deg b_{z,w}(\sq)\leq 0$ for all $z\in W$, where for a rational function $b(\sq)=P(\sq)/Q(\sq)$ we define $\deg b(\sq)=\deg P(\sq)-\deg Q(\sq)$. Fix $w\in W$ for the remainder of the proof. Let $z_0\in W$ be such that $\deg b_{z_0,w}(\sq)$ is maximal (taken over $z$ in the finite support of $\tw C_w$), and suppose for a contradiction that $\deg b_{z_0,w}(\sq)>0$. Define rational functions $c_{z,w}(\sq)$ by $\tw C_w=\sum_{z}c_{z,w}(\sq)T_z$. Then 
\begin{align}\label{eq:changebasis}
c_{z,w}(\sq)=b_{z,w}(\sq)+\sum_{z'\,:\,z'>z}b_{z',w}(\sq)P_{z,z'}(\sq^{-1}),
\end{align}
with $P_{z,z'}(\sq^{-1})$ the Kazhdan-Lusztig polynomials, and there are only finitely many nonzero terms in the sum. Since $\deg P_{z,z'}(\sq^{-1})\leq -1$ whenever $z'>z$ we have $\deg c_{z_0,w}(\sq)=\deg b_{z_0,w}(\sq)>0$. 

Now, by Lemma~\ref{lem:basic} we have $\tw^*=\tw$ and $\tw^2=1$ and so $\|\tw C_w\|^2=\|C_w\|^2$, and hence
\begin{align}\label{eq:norm2}
\|C_w\|^2=\sum_zc_{z,w}(\sq)^2. 
\end{align}
Since $\deg c_{z_0,w}(\sq)>0$ the right hand side is unbounded as $\sq\to\infty$. But on the other hand we have
$$
\|C_w\|^2=1+\sum_{z<w}P_{z,w}(\sq^{-1})^2,
$$
and hence $\|C_w\|^2\to 1$ as $\sq\to\infty$, a contradiction. 

Thus each rational function $b_{z,w}(\sq)$ is bounded as $\sq\to\infty$, and hence approaches a limit. Since $\tw$ is bar-invariant (by Lemma~\ref{lem:basic}) and $C_w$ is bar-invariant, each coefficient $b_{z,w}(\sq)$ is bar-invariant. Thus boundedness as $\sq\to\infty$ implies regularity at $\sq=0$, and moreover
\begin{align}\label{eq:limit1}
b_{z,w}(0)=\lim_{\sq\to\infty}b_{z,w}(\sq)\quad\text{for all $z\in W$}.
\end{align}
This shows that $
b_{z,w}(\sq)\in\ZZ[\sq,(1+\sq^2+\cdots+\sq^{2m-2})^{-1}\mid m\in M_I]
$. Moreover, taking limits in~(\ref{eq:changebasis}) gives 
$$
\lim_{\sq\to\infty}c_{z,w}(\sq)=\lim_{\sq\to\infty}b_{z,w}(\sq)=b_{z,w}(0),
$$
and so taking limits in~(\ref{eq:norm2}) we have
$$
1=\sum_{z\in W}b_{z,w}(0)^2. 
$$
Finally, Theorem~\ref{thm:actiononKLbasis} and Remark~\ref{rem:nokill} gives $b_{\rho_I^R(w),w}(\sq)=\eta_{I,R}(w)\in\{-1,1\}$, and thus all other terms $b_{z,w}(0)$ are zero, as required. 
\end{proof}

\subsection{The action of the cactus group on the asymptotic algebra}\label{sec:mainproof}

Recall the definition of Lusztig's homomorphism from Section~\ref{sec:LusHom}. On extending scalars, we consider $\psi$ as a homomorphism $\psi:\QQ(\sq)\otimes_{\sR} H\to \QQ(\sq)\otimes_{\ZZ}\cJ$. 

Let 
$$
\sR_I=\ZZ[\sq,(1+\sq^2+\cdots+\sq^{2m-2})^{-1}\mid m\in M_I]
$$
be the ring from Corollary~\ref{cor:coefficientring}. 

\begin{thm}\label{thm:mainresult}
Let $(W,S,\varphi)$ be a weighted Coxeter system with $S$ finite and $\varphi$ positive. Assume $\mathsf{P1}$--$\mathsf{P15}$, and assume that $\cD$ is finite. Let $I\subseteq S$ be spherical. Then:
\begin{compactenum}[$(1)$]
\item We have $\psi(\tw)\in\sR_I\otimes_{\ZZ}\cJ$. Thus $\psi(\tw)$ can be specialised at any real~$\sq$, and moreover the specialisations at $\sq=0$ and $\sq=1$ lie in the $\ZZ$-algebra~$\cJ$. 
\item Writing $\psi_0(\tw)\in\cJ$ for the specialisation of $\psi(\tw)$ at $\sq=0$, we have
\begin{compactenum}[$(a)$]
\item for every $w\in W$, 
$$
\psi_0(\tw)\st_w=\eta_{I,R}(w)\st_{\rho_I^R(w)}\quad\text{and}\quad \st_w\psi_0(\tw)=\eta_{I,L}(w)\st_{\la_I^L(w)}.
$$ 
\item the element $\psi_0(\tw)\in\cJ$ is given by
$$
\psi_0(\tw)=\sum_{d\in\cD}n_d\eta_{I,R}(d)\st_{\rho_I^R(d)}=\sum_{d\in\cD}n_d\eta_{I,L}(d)\st_{\la_I^L(d)}. 
$$
\item $\psi_0(\tw)^*=\psi_0(\tw)$. 
\end{compactenum}
\end{compactenum}
\end{thm}

\begin{proof}
With $b_{z,w}(\sq)$ as in~(\ref{eq:expand}), Lemma~\ref{lem:51} gives
\begin{align}\label{eq:unspec}
\psi(\tw)\st_w=\sum_{z\sim_{LR}w}b_{z,w}(\sq)\st_z. 
\end{align}
Since $1_{\cJ}=\sum_{d\in\cD}n_d\st_d$ it follows that 
$$
\psi(\tw)=\psi(\tw)1_{\cJ}=\sum_{d\in\cD}n_d\psi(\tw)\st_d=\sum_{z\in W}\bigg(\sum_{d\in\cD,\,d\sim_{LR}z}n_db_{z,d}(\sq)\bigg)\st_z.
$$
Thus by Theorem~\ref{thm:specialisation} the element $\psi(\tw)$ admits a specialisation at all real values of~$\sq$. From the description of the ring $\sR_I$ it is clear that the specialisation at $\sq=0$ lands in $\cJ$. For the specialisation at $\sq=1$, note that $\tw$ and $T_{\sw_I}$ both specialise to $\sw_I$, and so $\psi(\tw-T_{\sw_I})$ specialises to~$0$ at $\sq=1$. Thus $\psi(\tw)|_{\sq=1}=\psi(T_{\sw_I})|_{\sq=1}$, and the latter lies in $\cJ$ because Lusztig's homomorphism on $H$ is $\psi:H\to\ZZ[\sq,\sq^{-1}]\otimes_{\ZZ}\cJ$. 

In particular at $\sq=0$ Theorem~\ref{thm:specialisation} gives
$$
\psi_0(\tw)=\sum_{z\in W}\bigg(\sum_{d\in\cD,\,d\sim_{LR}z}n_db_{z,d}(0)\bigg)\st_z=\sum_{d\in \cD}n_d \eta_{I,R}(d)\st_{\rho_I^R(d)}
$$
(recall that Bonnaf\'e's right cactus action preserves left cells, and hence two-sided cells, and so 
$
\rho_I^R(d)\sim_{LR}d 
$). This proves the first equation in (2)(b). The first equation in (2)(a) is then immediate from~(\ref{eq:unspec}) and Theorem~\ref{thm:specialisation}. 

Recall the definition of the bilinear form $\langle\cdot,\cdot\rangle_{\infty}$ from Section~\ref{sec:KL}. Using the proved formula from (2)(a) we have $\langle\psi_0(\tw)\st_x,\psi_0(\tw)\st_y\rangle_{\infty}=\eta_{I,R}(x)\eta_{I,R}(y)\delta_{\rho_I^R(x),\rho_I^R(y)}=\delta_{x,y}$, and hence $\langle \psi_0(\tw)a,\psi_0(\tw)b\rangle_{\infty}=\langle a,b\rangle_{\infty}$ for all $a,b\in\cJ$. Thus $\langle \psi_0(\tw)^*\psi_0(\tw)a,b\rangle_{\infty}=\langle a,b\rangle_{\infty}$ for all $a,b\in\cJ$, and nondegeneracy gives $\psi_0(\tw)^*\psi_0(\tw)a=a$. Taking $a=1_{\cJ}$ we have $\psi_0(\tw)^*\psi_0(\tw)=1_{\cJ}$, and since $\psi_0(\tw)^2=\psi_0(\tw^2)=1_{\cJ}$ we have $\psi_0(\tw)^*=\psi_0(\tw)$, proving~(c). 

The ``right-handed'' versions of (2)(a) and (2)(b) now follow by applying $*$ to the ``left-handed'' equations, and using $\psi_0(\tw)^*=\psi_0(\tw)$ and the facts that $(\la_I^L)^{\mathsf{op}}=\rho_I^R$ and $\eta_{I,L}(w)=\eta_{I,R}(w^{-1})$ (see Theorem~\ref{thm:bonnafe}). 
\end{proof}

\begin{cor} \cite[Conjecture~4.2]{RW:24} holds. 
\end{cor}

\begin{proof}
In the equal parameter case the formulae in Theorem~\ref{thm:mainresult} reduce to those in~\cite[Conjecture~4.2]{RW:24}, because $n_d=1$, $\eta_I(x)=(-1)^{\ba(\sw_Ix)}$, $\rho_I=\sigma_{W_I}$, and $\lambda_I=\omega_I\circ\rho_I$. 
\end{proof}

As discussed in the introduction, and using Theorem~\ref{thm:mainresult}, there are group homomorphisms 
$$\theta_1:\cact(W)\to (\QQ(\sq)\otimes_{\sR}H)^{\times},\quad \theta_2:\cact(W)\to (\QQ(\sq)\otimes_{\ZZ}\cJ)^{\times},\quad\text{and}\quad \theta_3:\cact(W)\to \cJ^{\times}$$
with $\theta_1(\gamma_I)=\tw$, $\theta_2(\gamma_I)=\psi(\tw)$, and $\theta_3(\gamma_I)=\psi_0(\tw)$ (specialisation at $\sq=0$). Following from \cite[Remark~4.4]{RW:24}, it appears to not be known whether $\theta_1$ or $\theta_2$ are faithful. Since Lusztig's homomorphism is injective (see \cite[Proposition~18.12]{Lus:03}) and $\theta_2=\psi\circ\theta_1$, faithfulness of $\theta_1$ is equivalent to faithfulness of~$\theta_2$. We show in the proposition below that the homomorphism $\theta_3$ is typically not faithful.

\begin{prop}\label{prop:nonfaithful}
If $W$ has at least one irreducible component of rank at least~$2$, then the homomorphism $\theta_3:\cact(W)\to \cJ^{\times}$ is not faithful. 
\end{prop}

\begin{proof}
By assumption there exists $s,s'\in S$ with $ss'\neq s's$. Since $\sq^{\varphi(s)}T_sC_1=-C_1+\sq^{\varphi(s)}C_s$ and $\sq^{-\varphi(s)}T_sC_s=C_s$ we have $\eta_{\{s\}}(1)=-1$, $\eta_{\{s\}}(s)=1$, $\rho_{\{s\}}(1)=1$, and $\rho_{\{s\}}(s)=s$ (cf. Theorem~\ref{thm:ML}). Thus for $w\in W$ we have $\rho_{\{s\}}^R(w)=w$ and 
$$
\eta_{\{s\},R}(w)=\begin{cases}
1&\text{if $\ell(sw)=\ell(w)-1$}\\
-1&\text{if $\ell(sw)=\ell(w)+1$},
\end{cases}
$$
and so Theorem~\ref{thm:mainresult} gives
$$
\psi_0(\widetilde{\sw}_{\{s\}})=\sum_{d\in\cD,\,sd<d}n_d\st_{d}-\sum_{d\in\cD,\,sd>d}n_d\st_{d}.
$$
Lusztig's properties imply that $\st_d\st_{d'}=\delta_{d,d'}n_d\st_d$ for $d,d'\in\cD$, and it follows that $\psi_0(\widetilde{\sw}_{\{s\}})$ and $\psi_0(\widetilde{\sw}_{\{s'\}})$ commute, and hence $\theta_3((\gamma_{\{s\}}\gamma_{\{s'\}})^2)=1$. However we have $(\gamma_{\{s\}}\gamma_{\{s'\}})^2\neq 1$ in the cactus group, for otherwise applying the homomorphism $\cact(W)\to W$ with $\gamma_I\mapsto \sw_I$ we have $ss'=s's$.
\end{proof}

\subsection{The $\tilde{\sA}_2$ example}\label{sec:example}

We conclude by computing the right action of $\cact(W)$ (corresponding to left multiplication on the asymptotic algebra) for the $\tilde{\sA}_2$ Coxeter group 
$$
W=\langle s_0,s_1,s_2\mid s_0^2=s_1^2=s_2^2=(s_0s_1)^3=(s_0s_2)^3=(s_1s_2)^3=1\rangle.
$$
For $i,j\in\{0,1,2\}$ with $i\neq j$ let $\sw_{ij}=s_is_js_i$ be the longest element of $W_{ij}=\langle s_i,s_j\rangle$. 

The Hecke algebra of $W$ necessarily has equal parameters, and so by Remark~\ref{rem:etaequal} the sign map $\eta_{I,R}:W\to\{-1,1\}$ is given by
\begin{align}\label{eq:signs}
\eta_{I,R}(w)=\eta_I(x)=(-1)^{\ba(\sw_I x)}
\end{align}
where $w=xy$ with $x\in W_I$ and $y\in{^I}W$ (for $I\subseteq S=\{s_0,s_1,s_2\}$ spherical).

The cactus group $\cact(W)$ is generated by elements $\gamma_i=\gamma_{\{s_i\}}$ ($0\leq i\leq 2$) and $\gamma_{ij}=\gamma_{ji}=\gamma_{\{s_i,s_j\}}$ ($0\leq i< j\leq 2$) subject to the relations 
$$
\gamma_i^2=1,\quad 
\gamma_{ij}^2=1,\quad\gamma_{ij}\gamma_i\gamma_{ij}=\gamma_j.
$$
The right action of $\cact(W)$ on $W$ is given by 
$$
\gamma_I^R(w)=\rho_I^R(w)=\rho_I(x)y,
$$
where $w=xy$ with $x\in W_I$ and $y\in {^I}W$, and $\rho_I:W_I\to W_I$ is the involution given by the Mathas-Lusztig Theorem. As noted in the proof of Proposition~\ref{prop:nonfaithful} we have $\gamma_i^R(w)=w$ for all $w\in W$. Therefore the main task is to compute $\gamma_{ij}^R(w)$. To do so, we must compute $\rho_{ij}(x)$ for $x \in W_{ij}$. This is a direct calculation in a Coxeter group of type $\sA_2$, and we have (see \cite[Proposition~6.5]{Bon:16})
\begin{align}\label{eq:paraboliccalculation}
\begin{tabular}{l||l|l|l|l|l|l|}
$w$&$1$&$s_i$&$s_j$&$s_is_j$&$s_js_i$&$s_is_js_i$\\
\hline\hline
$\rho_{ij}(w)$&$1$&$s_js_i$&$s_is_j$&$s_j$&$s_i$&$s_is_js_i$\\
\hline
\end{tabular}
\end{align}

\begin{figure}[H]
\centering
\begin{tikzpicture}[scale=0.9]
  \path [fill=blue!50] (0,{2*0.866})--({-0.5*2},{4*0.866})--({0.5*2},{4*0.866})--(0,{2*0.866});
   \path [fill=blue!50] (1,0)--(3.5,0)--(3.5,{-3*0.866})--(2.5,{-3*0.866})--(1,0);
    \path [fill=blue!50] (-1,0)--(-3.5,0)--(-3.5,{-3*0.866})--(-2.5,{-3*0.866})--(-1,0);
       \path [fill=blue!50] (0.5,0.866)--(3.5,0.866)--(3.5,{4*0.866})--(2,{4*0.866})--(0.5,0.866);
        \path [fill=blue!50] (-0.5,0.866)--(-3.5,0.866)--(-3.5,{4*0.866})--(-2,{4*0.866})--(-0.5,0.866);
         \path [fill=blue!50] (0,0)--(1.5,{-3*0.866})--(-1.5,{-3*0.866})--(0,0);
         \path [fill=yellow!70] (0,{2*0.866})--(-1,{4*0.866})--(-2,{4*0.866})--(-0.5,{1*0.866})--(0.5,{1*0.866})--(2,{4*0.866})--(1,{4*0.866})--(0,{2*0.866});
         \path [fill=yellow!70] (0.5,{1*0.866})--(3.5,{1*0.866})--(3.5,{0*0.866})--(1,{0*0.866})--(2.5,{-3*0.866})--(1.5,{-3*0.866})--(0,{0*0.866})--(0.5,{1*0.866});
         \path [fill=yellow!70] (-0.5,{1*0.866})--(-3.5,{1*0.866})--(-3.5,{0*0.866})--(-1,{0*0.866})--(-2.5,{-3*0.866})--(-1.5,{-3*0.866})--(0,{0*0.866})--(-0.5,{1*0.866});
    \path [fill=black!60] (0,0) -- (-0.5,0.866) -- (0.5,0.866) -- (0,0);
    \draw (2.5, {-3*0.866})--( 3.5, {-1*0.866} );
    \draw (1.5, {-3*0.866})--( 3.5, {1*0.866} );
    \draw  (0.5, {-3*0.866})--( 3.5, {3*0.866} );
    \draw  (-0.5, {-3*0.866})--( 3, {4*0.866} );
    \draw  (-1.5, {-3*0.866})--( 2, {4*0.866} );
    \draw   (-2.5, {-3*0.866})--( 1, {4*0.866} );
    \draw  (-3.5, {-3*0.866})--(0, {4*0.866} );
    \draw  (-3.5, {-1*0.866})--(-1, {4*0.866} );
    \draw  (-3.5, {1*0.866})--(-2, {4*0.866} );
    \draw  (-3.5, {3*0.866})--(-3, {4*0.866} );
   \draw (-2.5, {-3*0.866})--( -3.5, {-1*0.866} );
    \draw (-1.5, {-3*0.866})--( -3.5, {1*0.866} );
    \draw  (-0.5, {-3*0.866})--( -3.5, {3*0.866} );
    \draw  (0.5, {-3*0.866})--( -3, {4*0.866} );
    \draw (1.5, {-3*0.866})--( -2, {4*0.866} );
    \draw  (2.5, {-3*0.866})--( -1, {4*0.866} );
    \draw  (3.5, {-3*0.866})--(0, {4*0.866} );
    \draw  (3.5, {-1*0.866})--(1, {4*0.866} );
    \draw  (3.5, {1*0.866})--(2, {4*0.866} );
    \draw  (3.5, {3*0.866})--(3, {4*0.866} );
    \draw (-3.5, -2.598)--( 3.5, -2.598);
    \draw (-3.5, -1.732)--( 3.5, -1.732);
    \draw (-3.5, -0.866)--( 3.5, -0.866);
    \draw (-3.5, 0)--( 3.5, 0);
    \draw (-3.5, 3.464)--( 3.5, 3.464 );
    \draw (-3.5, 2.598)--( 3.5, 2.598);
    \draw (-3.5, 1.732)--( 3.5, 1.732);
    \draw (-3.5, 0.866)--( 3.5, 0.866);
\end{tikzpicture}
\caption{Two-sided cells and right cells in $\tilde{\sA}_2$}\label{fig:A2partitions}
\end{figure}
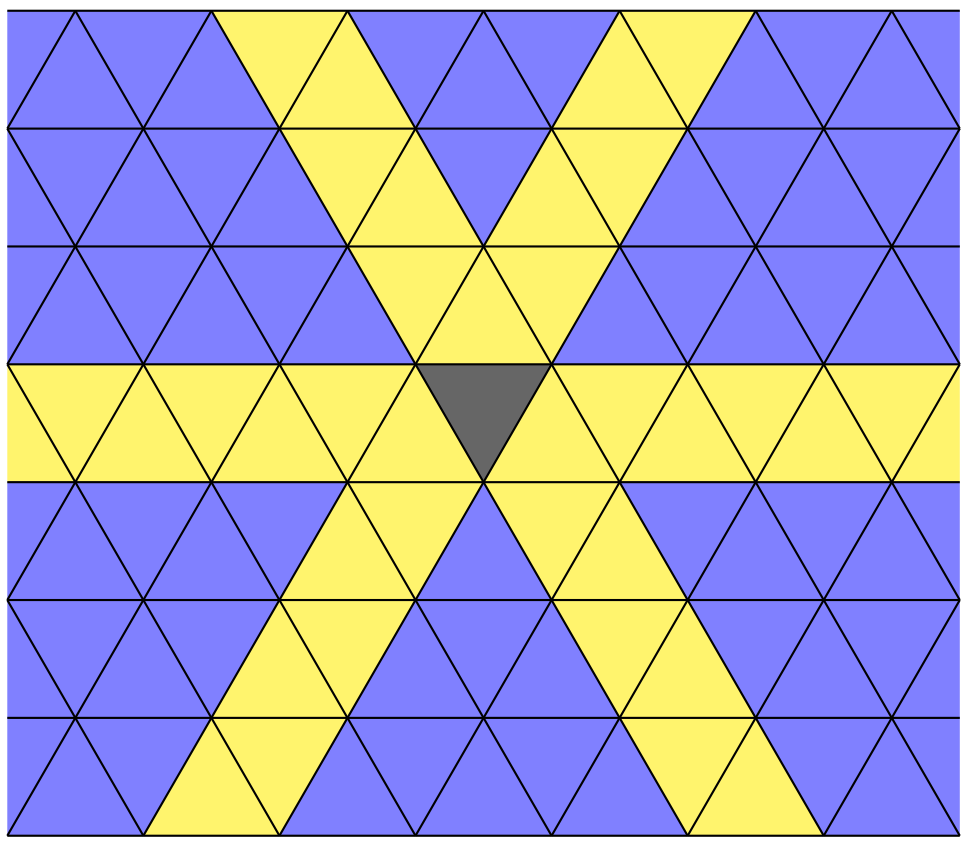

The decomposition of $W$ into cells is illustrated in Figure~\ref{fig:A2partitions} (see \cite{Lus:I} and \cite[Chapter~23]{Bon:17}). Elements of the same colour lie in the same two-sided cell, and connected components are right cells. Left cells are the inverses of right cells. Thus there are precisely $3$ two-sided cells: the \textit{lowest} two-sided cell $\Delta_0$ (consisting of the 6 blue right cells), the \textit{middle} two-sided cell $\Delta_1$ (consisting of the $3$ yellow right cells) and the \textit{highest} two-sided cell $\Delta_2=\{1\}$ (the grey alcove). The $\ba$-function is constant on two-sided cells, with values $\ba(\Delta_0)=3$, $\ba(\Delta_1)=1$, and $\ba(\Delta_2)=0$.

By inspection we have
$$
\Delta_0=\{\sw_{ij}y,s_k\sw_{ij}y\mid \{i,j,k\}=\{0,1,2\}\text{ and $y\in W$ with $D_L(y)\subseteq \{s_k\}$}\},
$$
where $D_L(y)$ is the left descent set of~$y$.

Let $\cR_0,\cR_1,\cR_2$ be the three right cells of $\Delta_1$, with $s_i\in\cR_i$ for $i=0,1,2$. Each cell $\cR_i$ is naturally in bijection with $\ZZ$, with $s_i\leftrightarrow 0$, and $s_is_{i+1}\leftrightarrow 1$ (where indices are read modulo~$3$). Thus the two-sided cell $\Delta_1$ is naturally in bijection with $\{0,1,2\}\times \ZZ$, and we write $x_{i,n}$ for the element of $\Delta_1$ corresponding to $(i,n)\in\{0,1,2\}\times\ZZ$. For example $x_{2,-4}=s_2s_1s_0s_2s_1$ and $x_{1,3}=s_1s_2s_0$. 

Since the cactus action preserves two-sided cells, we consider the action restricted to each of the two-sided cells. The action on the highest cell $\Delta_2$ is trivial, and so we focus on the lowest and middle two-sided cells in the propositions below.

\begin{prop}\label{prop:Delta0}
Let $\{i,j,k\}=\{0,1,2\}$ and let $y\in W$ with $D_L(y)\subseteq \{s_k\}$. The right action of the cactus group on $\Delta_0$ is given by the following table. 
\begin{center}
\begin{tabular}{l||l|l|l|}
$w$&$\gamma_{ij}^R(w)$&$\gamma_{ik}^R(w)$&$\gamma_{jk}^R(w)$\\
\hline\hline
$\sw_{ij}y$&$\sw_{ij}y$&$s_k\sw_{ij}y$&$s_k\sw_{ij}y$\\
\hline
$s_k\sw_{ij}y$&$s_k\sw_{ij}y$&$\sw_{ij}y$&$\sw_{ij}y$\\
\hline
\end{tabular}
\end{center}
\end{prop}

\begin{proof}
Consider $w=\sw_{ij}y$. We have $\pr_{ij}^R(w)=\sw_{ij}$, $\pr_{ik}^R(w)=s_i$, and $\pr_{jk}^R(w)=s_j$. Thus by~(\ref{eq:paraboliccalculation}) we have $\gamma_{ij}^R(w)=\rho_{ij}(\sw_{ij})y=\sw_{ij}y$, $\gamma_{ik}^R(w)=\rho_{ik}(s_i)s_js_iy=s_k\sw_{ij}y$, and $\gamma_{jk}^R(w)=\rho_{jk}(s_j)s_is_jy=s_k\sw_{ij}y$. The calculations for $w=s_k\sw_{ij}y$ are similar. 
\end{proof}

\begin{prop}\label{prop:Delta1}
Let $i\in\{0,1,2\}$, $j=i+1\pmod{3}$, $k=i+2\pmod{3}$, and $n\in\ZZ$. The right action of the cactus group on $\Delta_1$ is given by:
\begin{align*}
\gamma_{ij}^R(x_{i,n})=x_{j,n-1},\quad
\gamma_{ij}^R(x_{j,n})=x_{i,n+1},\quad\text{and}\quad 
\gamma_{ij}^R(x_{k,n})=x_{k,n}.
\end{align*}
\end{prop}

\begin{proof}
If $n\geq1$ then $x_{i,n}=s_is_js_ks_is_js_k\cdots$ (length $n+1\geq 2$) and hence by~(\ref{eq:paraboliccalculation}) we have
$$
\gamma_{ij}^R(x_{i,n})=\rho_{ij}(s_is_j)s_ks_is_js_k\cdots=(s_j)s_ks_is_js_k\cdots=x_{j,n-1}.
$$
If $n\leq -1$ then $x_{i,n}=s_is_ks_js_is_ks_j\cdots$ (length $|n|+1\geq 2$) and hence 
$$
\gamma_{ij}^R(x_{i,n})=\rho_{ij}(s_i)s_ks_js_is_ks_j\cdots=(s_js_i)s_ks_js_is_ks_j\cdots=x_{j,n-1}.
$$
If $n=0$ then $x_{i,0}=s_i$, and we have $\gamma_{ij}^R(x_{i,0})=\rho_{ij}(s_i)=s_js_i=x_{j,-1}$. Thus $\gamma_{ij}^R(x_{i,n})=x_{j,n-1}$ for all $n\in\ZZ$. The remaining calculations are similar. 
\end{proof}

Let $\mathsf{orb}(w)=\{\gamma^R(w) \mid \gamma\in\cact(W)\}$ be the orbit of $w\in W$ under the right cactus action. Let $\cR(w)$ (respectively $\cL(w)$) denote the right (respectively left) cell of $W$ containing~$w$. Recall that $\rho_I^R(w)\sim_L w$ for all $w\in W$, and hence for each $w\in W$ we have
\begin{align}\label{eq:orb}
\mathsf{orb}(w)\subseteq\cL(w).
\end{align}
In \cite[Theorem~1.1]{Los:19} it is shown that in finite type $\sA_n$ the orbits of the right action are precisely the left cells. That is, $\mathsf{orb}(w)=\cL(w)$. This phenomenon does not hold more generally, as illustrated in the corollary below. 

\begin{cor}
We have the following.
\begin{compactenum}[$(1)$]
\item If $w\in\Delta_0$ then $|\mathsf{orb}(w)|=2$. 
\item If $w\in \Delta_1$ then $\mathsf{orb}(w)=\cL(w)$. 
\end{compactenum}
Thus the right cactus group action is transitive on each left cell of~$\Delta_1$, and is not transitive on any left cell of $\Delta_0$.
\end{cor}

\begin{proof}
(1) follows from Proposition~\ref{prop:Delta0}. Indeed, each orbit is of the form $\{\sw_{ij}y,s_k\sw_{ij}y\}$ where $\{i,j,k\}=\{0,1,2\}$ and $D_L(y)\subseteq\{s_k\}$. In particular, since each left cell in $\Delta_0$ is infinite the action is not transitive on any left cell of $\Delta_0$. 

(2) Let $i\in\{0,1,2\}$, $j=i+1\pmod{3}$, and $k=i+2\pmod{3}$. We have 
$$
\cL(x_{i,n})=\cR(x_{i,n}^{-1})^{-1}=\{x_{i,n+3m},x_{j,n-1+3m},x_{k,n-2+3m}\mid m\in\ZZ\}. 
$$
It is straightforward to check, using Proposition~\ref{prop:Delta1}, that
$$
\cL(x_{i,n})=\{x_{i,n+3m},x_{j,n-1+3m},x_{k,n-2+3m}\mid m\in\ZZ\}\subseteq \mathsf{orb}(x_{i,n}),
$$
and so by~(\ref{eq:orb}) we have $\mathsf{orb}(x_{i,n})=\cL(x_{i,n})$. 
\end{proof}

Finally we illustrate Theorem~\ref{thm:mainresult} by computing $\psi_0(\tw)$ for $I=\{s_1,s_2\}$. Since $W$ has equal parameters we have $n_d=1$ for all $d\in\cD$. The elements $d\in\cD$, their images $\gamma_I^R(d)$, and signs $\eta_{I,R}(d)$ are given in the table below (using~(\ref{eq:signs}) and (\ref{eq:paraboliccalculation})). 
\begin{center}
\begin{tabular}{l||l|l|l|l|l|l|l|l|l|l|}
$d$&$1$&$s_0$&$s_1$&$s_2$&$\sw_{12}$&$\sw_{01}$&$\sw_{02}$&$s_0\sw_{12}s_0$&$s_2\sw_{01}s_2$&$s_1\sw_{02}s_1$\\
\hline
$\gamma_I^R(d)$&$1$&$s_0$&$s_2s_1$&$s_1s_2$&$\sw_{12}$&$s_2\sw_{01}$&$s_1\sw_{02}$&$s_0\sw_{12}s_0$&$\sw_{01}s_2$&$\sw_{02}s_1$\\
\hline
$\eta_{I,R}(d)$&$-1$&$-1$&$-1$&$-1$&$+1$&$-1$&$-1$&$-1$&$-1$&$-1$\\
\hline
\end{tabular}
\end{center}
Thus the formula in Theorem~\ref{thm:mainresult}(2)(b) gives
$$
\psi_0(\tw)=-\st_1-\st_{s_0}-\st_{s_2s_1}-\st_{s_1s_2}+\st_{\sw_{12}}-\st_{s_2\sw_{01}}-\st_{s_1\sw_{02}}-\st_{s_0\sw_{12}s_0}-\st_{\sw_{01}s_2}-\st_{\sw_{02}s_1}.
$$

\subsection*{Acknowledgements} The research of the second author is supported by ARC Discovery Project DP200100712. We thank Noah White for enlightening comments and suggestions, and the second author thanks Joel Kamnitzer for introducing him to the cactus group action on cells. OpenAI's Astra 6.0 was used to assist in reviewing and revising mathematical arguments and their presentation. The core strategy and proofs in the paper are due to the authors, and they take full responsibility for the mathematical correctness of the results presented. 

\bibliographystyle{plain}

\bigskip
  \footnotesize
  
   \noindent Fabricio Dos Santos, \textsc{School of Mathematics and Statistics, University of Sydney, NSW 2006, Australia}, \textit{E-mail address:} \texttt{f.dossantos@maths.usyd.edu.au}
   \smallskip
 
 \noindent James Parkinson, \textsc{School of Mathematics and Statistics, University of Sydney, NSW 2006, Australia}, 
 \noindent \textit{E-mail address:} \texttt{jamesp@maths.usyd.edu.au}

\end{document}